\documentclass{article}

\usepackage{amssymb,amsmath,color}
\usepackage{amsthm}
\usepackage{url}
\usepackage{tikz,subcaption}
\usepackage[enableskew]{youngtab}
\usepackage{ytableau,varwidth}
\usepackage{soul}

\definecolor{red}{rgb}{1,0,0}
\definecolor{blue}{rgb}{.2,.2,.8}

\def\mex{\mathrm{mex}}

\newtheorem{theorem}{Theorem}[section]
\newtheorem{corollary}[theorem]{Corollary}

\newtheorem{conjecture}{Conjecture}

\theoremstyle{definition}
\newtheorem{definition}{Definition}

\newtheorem{remark}{Remark}

\usepackage[round-mode=places,round-precision=0,group-separator={\,},output-decimal-marker={.}]{siunitx}

\begin{document}

 \title{$6$-regular partitions: new combinatorial properties, congruences, and linear inequalities }
	\author{Cristina Ballantine
	\\
	\footnotesize Department of Mathematics and Computer Science\\
	\footnotesize College of The Holy Cross\\
	\footnotesize Worcester, MA 01610, USA \\
	\footnotesize cballant@holycross.edu
	\and Mircea Merca
	\\ 
		\footnotesize Department of Mathematical Methods and Models\\
\footnotesize Fundamental Sciences Applied in Engineering Research Center\\ \footnotesize University Politehnica of Bucharest\\
\footnotesize RO-060042 Bucharest, Romania\\
\footnotesize mircea.merca@profinfo.edu.ro
}
	\date{}
	\maketitle


\begin{abstract}
We consider the number of the $6$-regular partitions of $n$, $b_6(n)$, and  give infinite families  of congruences modulo $3$ (in arithmetic progression) for $b_6(n)$.  We also consider the number of the partitions of $n$ into distinct parts not congruent to $\pm 2$  modulo $6$, $Q_2(n)$, and 
investigate connections between  $b_6(n)$ and $Q_2(n)$ providing new combinatorial interpretations for these partition functions. In this context, we discover new infinite families of linear inequalities involving Euler's partition function $p(n)$. Infinite families of linear inequalities involving the $6$-regular partition function $b_6(n)$ and the distinct partition function $Q_2(n)$ are proposed as open problems. 
\end{abstract}

\section{Introduction}

Recall that a partition of a positive integer $n$ is a sequence of positive integers whose sum is $n$. The order of the summands is unimportant when writing the partitions of $n$, but for consistency, a partition of $n$ will be written with the summands in a nonincreasing order \cite{Andrews98}. As usual, we denote by $p(n)$ the number of integer partitions of $n$ and we have the generating function
\begin{align*}
\sum_{n=0}^\infty p(n)\, q^n = \frac{1}{(q;q)_\infty}. 
\end{align*}
Here and throughout, we use the following customary $q$-series notation:
\begin{align*}
& (a;q)_n = \begin{cases}
1, & \text{for $n=0$,}\\
(1-a)(1-aq)\cdots(1-aq^{n-1}), &\text{for $n>0$;}
\end{cases}\\
& (a;q)_\infty = \lim_{n\to\infty} (a;q)_n.
\end{align*}
Moreover, we use the short notation
$$
(a_1,a_2,\ldots,a_n;q)_\infty = (a_1;q)_\infty (a_2;q)_\infty \cdots (a_n;q)_\infty.
$$
Because the infinite product $(a;q)_{\infty}$ diverges when $a\neq 0$ and $|q| \geqslant 1$, whenever
$(a;q)_{\infty}$ appears in a formula, we shall assume $|q| < 1$.

For an integer $\ell>1$,  a partition is called $\ell$-regular if none of its parts is divisible by $\ell$. The number of the $\ell$-regular partitions of $n$ is usually denoted by $b_\ell(n)$ and its  arithmetic propertys are investigated in many interesting papers 
by 
Z.~Ahmed and N.~D.~Baruah \cite{Ahmed},
R.~Carlson and J.~J.~Webb, \cite{Carlson}, 
S.-P.~Cui and N.~S.~S.~Gu \cite{Cui}, 
B.~Dandurand and D.~Penninston \cite{Dand},
D.~Furcy and D.~Penniston \cite{Furcy},
M.~D.~Hirschhorn and J.~A.~Sellers \cite{Hirschhorn},
Q.-H.~Hou, L.~H.~Sun and  L.~Zhang \cite{Hou},
J.~Lovejoy and D.~Penniston \cite{Lovejoy},
D.~Penniston \cite{Penn,Penn8},
E.~X.~W.~Xia \cite{Xia15},
E.~X.~W.~Xia and O.~X.~M. Yao \cite{Xia14},
L.~Wang \cite{Wang17,Wang},
and J.~J.~Webb \cite{Webb}.
Elementary techniques in the theory of partitions  give the generating function
\begin{align}
\sum_{n=0}^\infty b_\ell(n)\,q^n = \frac{(q^\ell;q^\ell)_\infty}{(q;q)_\infty}.\label{eq:1}
\end{align}

In $2010$, G.~E.~Andrews, M.~D.~Hirschhorn and J.~A.~Sellers \cite{Andrews10} proved that $b_4(n)$ satisfies two infinite families of congruences modulo $3$. 
After a year, J.~J.~Webb \cite{Webb} proved an analogous result for $b_{13}(n)$. 
In $2012$, D.~Furcy and D.~Penniston \cite{Furcy} extended these results to other values of $\ell$ which are congruent to $1$ modulo $3$, i.e., $\ell\in\{7,19,25,34,37,43,49\}$. All these congruences are of the form
$$
b_{\ell}(3^{\beta} n + d) \equiv 0 \pmod 3.
$$
In addition, D.~Furcy and D.~Penniston \cite{Furcy} proved that
$$
b_{10}(9n+3) \equiv b_{22}(27n+16) \equiv b_{28}(27n+9) \equiv 0 \pmod 3.
$$
More recently, in 2015, Q.-H.~Hou, L.~H.~Sun and  L.~Zhang \cite{Hou} found infinite families of congruence relations modulo $3$, $5$ and $7$ for $\ell$-regular partitions with $\ell\in \{3,5,6,7,10\}$. In particular, when $\ell=6$, they proved that for $\alpha, n$ nonnegative integers,  $p_i$ primes congruent to $13, 17, 19, 23\pmod{24}$ and $j\not \equiv 0 \pmod{p_{\alpha+1}}$, 
\begin{equation}\label{hou_6} 
b_6\left( p_1^2\cdots p_{\alpha+1}^2n+\frac{p_1^2\cdots p_\alpha^2p_{\alpha+1}(24j+5p_{\alpha+1})-5}{24}\right)\equiv 0 \pmod 3.
\end{equation} 
Then, setting  $\alpha=0$ in \eqref{hou_6}, it follows that  for all $n\geqslant0$, $p\equiv 13,17,19,23 \pmod {24}$ prime, and $j\not \equiv 0 \pmod{p}$,  
\begin{equation}\label{hou_6_0} 
b_6\left( p^2n+pj +5\,\frac{p^2-1}{24}\right)\equiv 0 \pmod 3.
\end{equation}  

It turns out that the result in \cite{Hou} can be extended to other choices of primes. 

\begin{theorem}\label{th1} Let $\alpha$ be a nonnegative integer and let  $p_i\geqslant 5$, $1\leqslant i \leqslant \alpha+1$  be primes. If   $p_{\alpha+1} \equiv 3\pmod 4$ and $j\not \equiv 0 \pmod{p_{\alpha+1}}$, then for all integers $n\geqslant 0$ we have
\begin{equation}\label{hou_6_gen} 
b_6\left( p_1^2\cdots p_{\alpha+1}^2n+\frac{p_1^2\cdots p_\alpha^2p_{\alpha+1}(24j+5p_{\alpha+1})-5}{24}\right)\equiv 0 \pmod 3.
\end{equation} 

\end{theorem}

In particular, if $\alpha =0$, Theorem \ref{th1} states that \eqref{hou_6_0} holds for all primes $p\equiv 3 \pmod 4$, $j\not\equiv 0\pmod p$ and $n\geqslant 0$. This statement can be reformulated as follows.

For a prime $p\geqslant 5$, we set $$\alpha_p:= 5\,\frac{p^2-1}{24} \mod p,$$ where by $a \mod m$ we mean the residue of $a$ modulo $m$. 
Equivalently, $$\alpha_p=\left\lfloor 5p^2/24\right\rfloor \mod p$$ and also $$\alpha_p=-5\cdot 24_p^{-1} \mod p,$$ where $24_p^{-1}$ is the inverse of $24$ modulo $p$. 

Then, from Theorem \ref{th1} with $\alpha=0$ and \eqref{hou_6_0}, we obtain the following result, 
\begin{corollary} If $p$ is a prime congruent to $7,11,13,17,19,23$  modulo $24$  and $0\leqslant j\leqslant p-1$, $j\neq \lfloor 5p/24 \rfloor $, then for all $n\geqslant 0$ we have 
\begin{equation*}b_6\left( p^2n+pj +\alpha_p \right)\equiv 0 \pmod 3.\end{equation*}
\end{corollary}



We  also consider the partitions of $n$ into distinct parts not congruent to $\pm 2$ modulo $6$  in order to provide other properties for the number of $6$-regular partitions of $n$. 

\begin{definition}
	Let $n$ be a nonnegative integer. We define:
	\begin{enumerate}
		\item [i)] $b_{6,e}(n)$ to be the number of $6$-regular partitions of $n$ into an even number of parts;
		\item [ii)] $b_{6,o}(n)$ to be the number of $6$-regular partitions of $n$ into an odd number of parts.
	\end{enumerate} 
\end{definition} 

Clearly $b_6(n)=b_{6,e}(n)+b_{6,o}(n)$. For example, the partitions of $7$ into parts that are not multiples of $6$ are:  
\begin{align*}
& (7),\ (5,2),\ (5,1,1),\ (4,3),\ (4,2,1), (4,1,1,1),\ (3,3,1),\ (3,2,2),\ (3,2,1,1),\\
& (3,1,1,1,1),\ (2,2,2,1),\ (2,2,1,1,1),\ (2,1,1,1,1,1),\ (1,1,1,1,1,1,1).
\end{align*}
We see that $b_6(7)=14$, $b_{6,e}(7)=6$ and $b_{6,o}(7)=8$.

\begin{definition}
	Let $n$ be a nonnegative integer. We define $Q_2(n)$ to be the number of partitions of $n$ into distinct parts which are not congruent to $\pm 2$ modulo $6$.
\end{definition}

For example, the partitions of $14$ into distinct parts not congruent to $\pm 2$ modulo $6$ are: 
\begin{align*}
& (13,1),\ (11,3),\ (10,3,1),\ (9,5),\ (8,5,1).  
\end{align*}
Thus, $Q_2(14)=5$. 
The standard methods for producing partition generating functions (cf. \cite[Ch. 1]{Andrews98}) reveal directly that
\begin{align}
\sum_{n=0}^\infty Q_2(n)\, q^n = (-q,-q^3,-q^5,-q^6;q^6)_\infty \label{eq:2}
\end{align}
and the expansion starts as
$$
1+q+q^3+q^4+q^5+2q^6+2q^7+2q^8+3q^9+3q^{10}+3q^{11}+5q^{12}+5q^{13}+5q^{14}+\cdots.
$$
We remark that the sequences $Q_2(n)$ is known and can be seen in the
On-Line Encyclopedia of Integer Sequence \cite[A328796]{Sloane}.

  The following result introduces a new combinatorial interpretation for the partition function $Q_2(n)$.

 \begin{theorem}\label{Th:1}
 	For $n\geqslant 0$, $(-1)^n\, Q_2(n) =  b_{6,e}(n)-b_{6,o}(n) $.
 \end{theorem}

As a corollary of this theorem, we deduce the following parity result.

\begin{corollary}
	For $n\geqslant 0$, $Q_2(n)$ and $b_{6}(n)$ have the same parity.
\end{corollary}
 
 In order to obtain other combinatorial interpretations for the $6$-regular partitions of $n$ and the partitions of $n$ into distinct parts not congruent to $\pm 2$ modulo $6$, we consider the following restricted partition functions.
 
\begin{definition}
	Let $n$ be a nonnegative integer. We define 
	\begin{enumerate}
		\item[i)] $c(n)$ to be the number of partitions of $n$ into parts which are not congruent to
		$0$, $\pm 2$, $\pm 20$, $\pm 22$, $24$ modulo $48$;
		\item[ii)] $d(n)$ to be the number of partitions of $n$ into parts which are not congruent to
		$0$, $\pm 4$, $\pm 10$, $\pm 14$, $24$ modulo $48$.
	\end{enumerate}
\end{definition}
 
 We have the following result.
 
\begin{theorem}\label{Th:3}
	Let $n$ be a nonnegative integer. Then
	\begin{enumerate}
		\item[(i)]  $Q_2(n) =  c(n)-d(n-2) $;
		\item[(ii)] $b_6(n) =  c(n)+d(n-2) $.
	\end{enumerate}
\end{theorem}

The following corollary is a consequence of Theorems \ref{Th:1} and \ref{Th:3}. This result introduces new combinatorial interpretations for the $6$-regular partition functions $b_{6,e}(n)$ and $b_{6,o}(n)$.

\begin{corollary}\label{Cor 1.5}
	For $n\geqslant 0$, 
	\begin{enumerate}
		\item[(i)] $b_{6,e}(n) = \begin{cases}
		c(n), & \text{if $n$ is even}\\
		d(n-2), & \text{if $n$ is odd;}
		\end{cases}$
		\item[(ii)] $b_{6,o}(n) = \begin{cases}
		c(n), & \text{if $n$ is odd}\\
		d(n-2), & \text{if $n$ is even.}
		\end{cases}$
	\end{enumerate}
\end{corollary}

From Corollay \ref{Cor 1.5} we can obtain other combinatorial interpretations for the restricted partition functions $c(n)$ and $d(n)$. 

\begin{definition}
	Let $n$ be a nonnegative integer. We define 
	\begin{enumerate}
		\item[i)] $b_{6,ee}(n)$ to be the number of $6$-regular partitions of $n$ with an even number of even parts;
		\item[ii)] $b_{6,eo}(n)$ to be the number of $6$-regular partitions of $n$ with an odd number of even parts.
	\end{enumerate}
\end{definition}

Clearly $b_6(n)=b_{6,ee}(n)+b_{6,eo}(n)$. For example, the $6$-regular partitions of $7$ with an even number of even parts are:  
\begin{align*}
& (7),\ (5,1,1),\ (4,2,1),\ (3,3,1),\ (3,2,2),\\
& (3,1,1,1,1),\ (2,2,1,1,1),\ (1,1,1,1,1,1,1),
\end{align*}
while the $6$-regular partitions of $7$ with an odd number of even parts are:
\begin{align*}
& (5,2),\ (4,3),\ (4,1,1,1),\ (3,2,1,1), (2,2,2,1),\ (2,1,1,1,1,1).
\end{align*}
We see that $b_{6,ee}(7)=8$ and $b_{6,eo}(7)=6$.

Since the parity of the number of odd parts in a partition of $n$ is determined by the parity of $n$, we have $b_{6,ee}(n)$ equals $b_{6,e}(n)$ (respectively $b_{6,o}(n)$) if $n$ is even (respectively odd); and similarly for $b_{6,eo}(n)$. Thus, we have the following equivalent form of 
 Corollary \ref{Cor 1.5}. 
\begin{corollary}
For $n\geqslant 0$
\begin{enumerate}
    \item [(i)] $c(n)=b_{6,ee}(n)$;
    \item [(ii)] $d(n)=b_{6,eo}(n+2)$.
\end{enumerate}
\end{corollary}

In \cite{Andrews12}, while investigating the truncated form of Euler's pentagonal number theorem, 
\begin{align}
(q;q)_\infty = \sum_{n=-\infty}^\infty (-1)^n\, q^{n(3n-1)/2},\label{eq:3}
\end{align}
G. E. Andrews and M. Merca introduced the partition function $M_k(n)$, which counts the number of partitions of $n$ where $k$ is the least positive integer that is not a part and there are more parts $>k$ than there are parts $<k$.  For instance, we have $M_3(18)=3$ because the three partitions in question are 
$$(5,5,5,2,1),\ (6,5,4,2,1),\ (7,4,4,2,1).$$ 
Recently, Xia and Zhao \cite{Xia} defined $\widetilde{P}_k(n)$ to be the number of partitions of $n$ 
in which every part $\leqslant k$ appears at least once
and the first part larger that $k$ appears at least $k + 1$ times. For example, $\widetilde{P}_2(17)=9$, and the partitions in question are 
\begin{align*}
& (5,3,3,3,2,1),\ (4,4,4,2,2,1),\ (4,4,4,2,1,1,1),\\
& (4,3,3,3,2,1,1),\ (3,3,3,3,2,2,1),\ (3,3,3,3,2,1,1,1),\\
& (3,3,3,2,2,2,1,1),\ (3,3,3,2,2,1,1,1,1),\ (3,3,3,2,1,1,1,1,1,1).
\end{align*}

Considering \eqref{eq:1}, we easily deduce that the $6$-regular partition function $b_6(n)$ is closely related to Euler's partition function $p(n)$, i.e.,
\begin{align}
b_6(n) = \sum_{j=-\infty}^\infty (-1)^j\,  p\big(n-3j(3j-1)\big).\label{eq:4}
\end{align}  
There are two more general results for which  identity \eqref{eq:4} is the limiting cases $k\to\infty$.

\begin{theorem}\label{Th:5}
	For $n\geqslant 0$, $k>0$,
	\begin{align*}
	(-1)^{k} \left( b_6(n) - \sum_{j=-(k-1)}^k (-1)^j\, p\big(n-3j(3j-1)\big) \right)
	= \sum_{j=0}^{\lfloor n/6 \rfloor} b_6(n-6j)\, M_k(j). 
	\end{align*}
\end{theorem}

\begin{theorem}\label{Th:5a}
	For $n\geqslant 0$, $k>0$,
	\begin{align*}
	(-1)^{k-1} \left(b_6(n) - \sum_{j=-k}^k (-1)^j\, p\big(n-3j(3j-1)\big) \right)
	= \sum_{j=0}^{\lfloor n/6 \rfloor} b_6(n-6j)\, \widetilde{P}_k(j). 
	\end{align*}
\end{theorem}

On the other hand, by \eqref{eq:1} and \eqref{eq:3}, we can easily derive a linear recurrence relation similar to the Euler recurrence relation for $p(n)$, i.e.,
\begin{align}
\sum_{j=-\infty}^\infty (-1)^j\, b_6\big(n-j(3j-1)/2\big)
= \begin{cases}
(-1)^k, & \text{if $n=3k(3k-1)$, $k\in\mathbb{Z}$,}\\
0, & \text{otherwise.}
\end{cases}\label{eq:5}
\end{align}
 \begin{remark}
If we denote by $p_{<6}(n)$  the number of partitions of $n$ in which parts occur at most five times, Glaischer's bijection shows combinatorially that $b_6(n)=p_{<6}(n)$,  for all $n\geqslant 0$.   Then, the first proof of Theorem 1.1 in \cite{BW} with $4$ replaced by $6$ gives a combinatorial proof of \eqref{eq:5}. 
\end{remark}

Apart from this recurrence relation,
there is another linear recurrence relation for  $b_6(n)$.
For any integer $k$, let 
\begin{align*}
\rho_k 
:= \begin{cases}
-2, & \text{if $k\equiv 1 \pmod 3$,}\\
1, & \text{otherwise.}
\end{cases}
\end{align*}

\begin{theorem}\label{Th:6}
	For $n\geqslant 0$,
	$$
	\sum_{j=0}^\infty \rho_j\, b_6\big(n-j(j+1)/2\big)
	= \begin{cases}
	(-1)^k, & \text{if $n=k(3k-2)$, $k\in\mathbb{Z}$,}\\
	0, & \text{otherwise.}
	\end{cases}
	$$
\end{theorem}

As a consequence of Theorem \ref{Th:6}, we remark the following parity result which involves the  generalized octagonal numbers, $n(3n\pm 2)$.

\begin{corollary}
	For $n\geqslant 0$,
	$$
	\sum_{j=-\infty}^\infty b_6\big(n-3j(3j-1)/2\big) \equiv 1 \pmod 2
	$$
	if and only if $n$ is a generalized octagonal number.
\end{corollary}

In analogy with \eqref{eq:4}, we have the following result which shows that the
partition function $Q_2(n)$ can be express in terms of Euler's partition function $p(n)$ in two different ways.
  
\begin{theorem}\label{Th:8}
	For $n\geqslant 0$,
	\begin{enumerate}
		\item[(i)] $\displaystyle{Q_2(n) = \sum_{j=0}^\infty \rho_j\, p\big(n-j(j+1) \big)}$;
		\item[(ii)] $\displaystyle{Q_2(n) = \sum_{j=-\infty}^{\infty} p\left(\frac{n-j(3j-2)}{3}\right) }$,\\
		\item[] where $p(x)=0$ when $x$ is not a nonnegative integer.
	\end{enumerate}
\end{theorem}

\begin{remark}Using the notation of Andrews and Newman \cite{AN}, $\mex_{2,2}(\lambda)$ denotes the smallest even positive integer that is not a part of $\lambda$. 
We denote by $pm_{2j}(n)$ (respectively $pm_{>2j}(n)$) the number of partitions $\lambda$ of $n$ with $\mex_{2,2}(\lambda)=2j$ (respectively $\mex_{2,2}(\lambda)>2j$).  If $\lambda$ is a partition of $n-j(j+1)$ then $\lambda\cup(2j, 2(j-1), \ldots, 4,2)$ is a partition of $n$ with $\mex_{2,2}(\lambda)>2j$. Hence $$p\big(n-j(j+1)\big)-p\big(n-(j+1)(j+2)\big)=pm_{2j}(n).$$  Then, Theorem \ref{Th:8} (i) is equivalent to the statement that $Q_2(n)$ equals the number of partitions $\lambda$ of  $n$ with $\mex_{2,2}(\lambda)\equiv 2 \pmod 6$ minus the number of partitions $\lambda$ of $n$ with $\mex_{2,2}(\lambda)\equiv 4 \pmod 6$.
\end{remark}

Theorem \ref{Th:8} (i)  allows us to derive the following congruence identities.

\begin{corollary}
	For $n\geqslant 0$, 
 \begin{enumerate}
      \item[(i)] $\displaystyle{ \sum_{j=-\infty}^\infty p\big(n-3j(3j-1) \big)  \equiv Q_2(n) \pmod 2}$;
	 \item[(ii)] $\displaystyle{ \sum_{j=0}^\infty \, p\big(n-j(j+1) \big) \equiv Q_2(n)  \pmod 3}$.
  \end{enumerate}
\end{corollary}

Theorem \ref{Th:8} (ii) can be considered an identity of Watson type.
More details about identities of Watson type can be found in \cite{BM}.

In analogy with Theorem \ref{Th:6}, we have the following linear recurrence relations for the partition function $Q_2(n)$.

\begin{theorem}\label{Th:10}
	For $n\geqslant 0$,
	\begin{enumerate}
		\item[(i)] $\displaystyle{\sum_{j=-\infty}^\infty (-1)^j\, Q_2\big(n-j(3j-1)/2 \big) 
		= \begin{cases}
		\rho_k,& \text{if $n=k(k+1)$, $k\in\mathbb{N}_0$ }\\
		0,& \text{otherwise;}
		\end{cases}}$
		\item[(ii)] $\displaystyle{\sum_{j=-\infty}^\infty (-1)^j\, Q_2\big(n-3j(3j-1)/2 \big) 
	= \begin{cases}
	1,& \text{if $n=k(3k-2)$, $k\in\mathbb{Z}$ }\\
	0,& \text{otherwise.}
	\end{cases}}$
	\end{enumerate}
\end{theorem}

Theorem \ref{Th:10} (ii) provides a simple and
reasonably efficient way to compute the value of $Q_2(n)$. 
The number of terms in this linear recurrence relation is about $\sqrt{8n/9}$.
In fact, computing the value of $Q_2(n)$ with this linear recurrence relation requires all the values of $Q_2(k)$ with $k < n$.

The rest of this paper is organized as follows. 
Theorem \ref{th1} will be proved in Section \ref{S2a}.
In Sections \ref{S2}-\ref{S6}, we will provide
proofs of Theorems \ref{Th:1}, \ref{Th:3}, \ref{Th:5}, \ref{Th:6}, \ref{Th:8} and \ref{Th:10}. Our proof of these theorems rely on generating functions. For Theorems \ref{Th:1}, \ref{Th:8} (ii), and \ref{Th:10} (ii) we also give combinatorial proofs. (It would be very interesting to find  combinatorial proofs for the remaining theorems.) In the last section of this paper, we  propose as conjectures  new infinite families of linear inequalities for the partition functions $p(n)$, $b_6(n)$, and $Q_2(n)$.

\section{Proof of Theorem \ref{th1}}
\label{S2a}

Let 
$$f(-q)=\sum_{n=-\infty}^\infty(-1)^nq^{\frac{n(3n+1)}{2}}=(q;q)_\infty$$ 
and 
$$\psi(q)=\sum_{n=0}^\infty q^{\binom{n+1}{2}}=\frac{(q^2;q^2)_\infty}{(q;q^2)_\infty}$$ 
be Ramanujan's theta functions.  As mentioned in \cite{Hou} and also easily seen directly,  
\begin{align*}
\sum_{n=0}^{\infty}b_6(n)q^n& =\prod_{n=1}^{\infty}\frac{(q^6;q^6)_\infty}{(q;q)_\infty}\\ 
& \equiv \frac{(q^2;q^2)^2_\infty}{(q;q^2)_\infty}\pmod 3\\ 
& = f(-q^2)\psi(q) \pmod 3. 
\end{align*}
We rewrite the above expression as 
$$f(-q^2)\psi(q)= (q^2;q^2)_\infty \frac{(q^2;q^4)_\infty(q^4;q^4)_\infty}{(q;q^2)_\infty}=(q^2;q^2)_\infty(q^4;q^4)_\infty(-q;q^2)_\infty $$ 
Replacing $q$ by $-q$, we obtain 
$$f(-q^2)\psi(-q)= \frac{(q^2;q^2)^2_\infty}{(-q;q^2)_\infty}= (q^2;q^2)_\infty(q^4;q^4)_\infty(q;q^2)_\infty=(q;q)_\infty(q^4;q^4)_\infty. $$ 
Hence, $$\sum_{n=0}^{\infty}b_6(n)(-q)^n \equiv (q;q)_\infty(q^4;q^4)_\infty \pmod 3.$$

Let $\alpha$ be a nonnegative integer. Suppose $p_i\geqslant 5$, $1\leqslant i \leqslant \alpha+1$  are primes,   $p_{\alpha+1} \equiv 3\pmod 4$ and $j\not \equiv 0 \pmod{p_{\alpha+1}}$.  Given $n\geqslant 0$, we set $$m_n:=p_1^2\cdots p_{\alpha+1}^2n+\frac{p_1^2\cdots p_\alpha^2p_{\alpha+1}(24j+5p_{\alpha+1})-5}{24}.$$
 We show that for all nonnegative integers $n$ the coefficient of $q^{m_n}$ in $f(-q^2)\psi(-q)$ is zero.
We use Euler's pentagonal number theorem \eqref{eq:3} twice to see that 
$$f(-q^2)\psi(-q)=(q;q)_\infty(q^4;q^4)_\infty=\sum_{i,j=-\infty}^\infty (-1)^{i+j}q^{\frac{i(3i+1)}{2}+4\cdot \frac{j(3j+1)}{2}}.$$
We consider the equation $\frac{i(3i+1)}{2}+4\cdot \frac{j(3j+1)}{2}= m_n$ which is equivalent to \begin{align}\label{eq-dio} a^2+(2b)^2=24m_n+5\end{align} with $a=6i+1$, $b=6j+1$.

Since $j\not \equiv 0 \pmod{p_{\alpha+1}}$, it follows that $p_{\alpha+1}\mid 24m_n+5$ and $p_{\alpha+1}$ appears in the factorization of $24m_n+5$ with odd exponent. Then,  equation \eqref{eq-dio} has no solution and the coefficient of $q^{m_n}$ in $f(-q^2)\psi(-q)$ is zero. Hence $b_6(m_n)\equiv 0 \pmod 3$.  This concludes the proof of Theorem \ref{th1}.

\begin{remark} In the poof of Theorem \ref{th1} we reduced the congruence problem to a question of  representing  an integer as a sum of two squares. The proof of \cite[Theorem 2.3]{Hou} relies on a different Diophantine equation.  We note that    \cite[Theorem 2.2]{Hou} is a congruence result for  $b_3(n)$, the number of $3$-regular partitions of $n$. It is easy to see that the proof of \cite{Hou} reduces to representing an integer as the sum of two squares. Thus the same argument as in the proof of Theorem \ref{th1} can be used to show that the congruence modulo $3$ in \cite[Theorem 2.2]{Hou} holds in greater generality, i.e., only the prime $p_{\alpha+1}$ must be congruent to $3$ modulo $4$. For the convenience of the reader, we give the general statement below. 
\begin{theorem}
Let $\alpha$ be a nonnegative integer and let  $p_i\geqslant 5$, $1\leqslant i \leqslant \alpha+1$  be primes. If   $p_{\alpha+1} \equiv 3\pmod 4$ and $j\not \equiv 0 \pmod{p_{\alpha+1}}$, then for all integers $n\geqslant 0$ \begin{equation*} 
b_3\left( p_1^2\cdots p_\alpha^2p_{\alpha+1}^2n+\frac{p_1^2\cdots  p_\alpha^2p_{\alpha+1}(12j+p_{\alpha+1})-1}{12}\right)\equiv 0 \pmod 3.
\end{equation*} 
\end{theorem}
\end{remark}

\section{Proof of Theorem \ref{Th:1}}
\label{S2}
\subsection{Analytic proof}

Define
\begin{align*}
F(z,q) = \prod_{k=0}^\infty \frac{1}{(1-zq^{6k+1})(1-zq^{6k+2})(1-zq^{6k+3})(1-zq^{6k+4})(1-zq^{6k+5})}.
\end{align*}
On the other hand, we have
$$
F(z,q) = \sum_{m=0}^\infty \sum_{n=0}^\infty b_6(n,m)\,z^m\,q^n,
$$
where $b_6(n,m)$ is the number of partitions of $n$ with $m$ parts all of which are not congruent to $0$ modulo $6$.
 
Thus, considering the generating functions of $b_6(n,m)$ and $Q_2(n)$, we can write
\allowdisplaybreaks{
	\begin{align*}
F(-1,q) &= \sum_{n=0}^\infty \big( b_{6,e}(n) - b_{6,o}(n) \big)\, q^n
\end{align*}
and
\begin{align*}
F(-1,q) &=\frac{1}{(-q,-q^2,-q^3,-q^4,-q^5;q^6)_\infty} \\
& = \frac{(-q^6;q^6)_\infty}{(-q;q)_\infty} \\
& = (q;q^2)_\infty\,(-q^6;q^6)_\infty\\
& = \sum_{n=0}^\infty (-1)^n\, Q_2(n)\, q^n,
\end{align*}	
	where we have invoked the Euler identity \cite[(1.2.5)]{Andrews98}
	$$\frac{1}{(q;q^2)_\infty}=(-q;q)_\infty.$$

\subsection{Combinatorial proof} 
We remark first that the set $(S_1, S_2)$ with $S_1=\{n\in \mathbb N : n\not \equiv 0\pmod 6\}$ and $S-2=\{n\in \mathbb N : n\not \equiv \pm 2\pmod 6\}$ is not an Euler pair and the statement of Theorem \ref{Th:1} is not a special case of Theorem  3.1 of \cite{BMR}. However, the ideas used in the proof of \cite[Theorem 3.1]{BMR} can be used here. 
Given a partition $\lambda$, denote by $\ell(\lambda)$ the number of parts in $\lambda$. 
Note that in a $6$-regular partition, even parts are congruent to $\pm 2$ modulo $6$. 

Let $\mathcal B'_{6}(n)$ be the set of 6-regular partitions $\lambda$ of $n$  such that $\lambda$ has  at least one even part or  at least one repeated part which is not congruent to $3$ modulo $6$.  Moreover, denote by $\mathcal B'_{6,e}(n)$, respectively $\mathcal B'_{6,o}(n)$, the subset of partitions in $\mathcal B'_{6}(n)$ with $\ell(\lambda)$ even, respectively odd. We define an involution $\varphi$ on $\mathcal B'_{6}(n)$ that reverses the parity of $\ell(\lambda)$.

Start with $\lambda\in \mathcal B'_{6}(n)$. We denote by $r$ the largest repeated part of $\lambda$ that is not congruent to $3$ modulo $6$ and by $e$ the largest even part of $\lambda$.  If $r$ or $e$ do not exist, we set them equal to $0$. 
\begin{enumerate}
\item If $2r>e$,  we define $\varphi(\lambda)$ to be the partition obtained from $\lambda$ by replacing two parts equal to $r$ by a single part equal to $2r$. Note that, since $r\not \equiv 3 \pmod 6$, we have  $2r\not \equiv 0 \pmod 6$. Thus, $\varphi(\lambda)\in \mathcal B'_{6}(n)$.

\item If $2r\leqslant e$, we define $\varphi(\lambda)$ to be the partition obtained from $\lambda$ by replacing one part equal to $e$ by  two  parts equal to $e/2$. Note that since $e \equiv \pm 2\pmod 6$, we have $e/2\not \equiv 0,3 \pmod 6$. Thus, $\varphi(\lambda)\in \mathcal B'_{6}(n)$.
\end{enumerate}

Since $\varphi: \mathcal B'_{6}(n)\to \mathcal B'_{6}(n)$ is an involution that reverses the parity of $\ell(\lambda)$,  we have  that $|\mathcal B'_{6,e}(n)|=|\mathcal B'_{6,o}(n)|$. 

Let $\mathcal Q'_2(n)$ be the set of partitions $\lambda\in \mathcal B_6(n)$ with odd parts  and only parts congruent to $3$ modulo $6$ may be repeated. Since all parts of $\lambda$ are odd, $\ell(\lambda)\equiv n \pmod 2$. Thus, $$b_{6,e}(n)-b_{6,o}(n)=(-1)^n|\mathcal Q'_2(n)|.$$

Finally, we create a bijection $\psi:\mathcal Q'_2(n) \to \mathcal Q_2(n)$, where  $\mathcal Q_2(n)$ is the set of partitions of $n$ with distinct parts not congruent to $\pm 2$ modulo $6$. Here and throughout, if $k$ is a positive integer and  $\eta$ is a partition with  all parts divisible by $k$, we write $\eta_{/k}$ for the partition whose parts are the parts of $\eta$ divided by $k$. For any partition $\eta$, we denote by $k\eta$ the partition whose parts are the parts of $\eta$ multiplied by $k$. 

Let $\lambda \in \mathcal Q'_2(n)$. Here and throughout, by the union of two partition we mean the union of their multisets of parts arranged in nondecreasing order.   Write $\lambda=(\alpha, \beta)$ where $\alpha \cup \beta=\lambda$, $\alpha$ is a partition into distinct parts, and $\beta$ is a partition whose parts have even multiplicity. Thus, all parts of $\beta$ are congruent to $3$ modulo $6$ and  
 $\beta_{/3}$ is a partition into odd parts each with even multiplicity. We denote by  $\varphi_{Gl}$  Glaisher's bijection which maps a partition of $n$ with odd parts to a partition of $n$ into distinct parts. 
Then, $\varphi_{Gl}(\beta_{/3})$ is a partition  with even distinct parts. The partition $3\varphi_{Gl}(\beta_{/3})$ has distinct parts all congruent to $0$ modulo $6$. Set $\psi(\lambda): =\alpha \cup 3\varphi_{Gl}(\beta_{/3})$. Then, $\psi$ is a bijection from  $\mathcal Q'_2(n)$ to $\mathcal Q_2(n)$, which completes the proof of the theorem.  


\begin{remark} We note that, in fact,  the involution $\varphi$ in the combinatorial proof above reverses the parity of the number of even parts of a partition. Hence, the combinatorial proof above is also a proof for the following corrolary of  Theorem \ref{Th:1}. 

\begin{corollary}  For $n\geqslant 0$, $Q_2(n) =  b_{6,ee}(n)-b_{6,eo}(n) $.\end{corollary}
\end{remark}

\section{Proof of Theorem \ref{Th:3}}
\label{S3}

The Watson quintuple product identity \cite{Carlitz,Subbarao} states that
\begin{align}
\sum_{n=-\infty}^{\infty} q^{n(3n+1)/2}\,(z^{-3n}-z^{3n+1})
=(z,q/z,q;q)_{\infty}\,(qz^2,q/z^2;q^2)_{\infty}.\label{eq:8}
\end{align}
Elementary techniques in the theory of partitions give the following generating
functions 
\begin{align*}
 \sum_{n=0}^\infty c(n)\, q^n 
& = \frac{(q^2,q^{20},q^{22},q^{24},q^{26},q^{28},q^{46},q^{48};q^{48})_\infty}{(q;q)_\infty}\\
& = \frac{(q^2,q^{22},q^{24};q^{24})_\infty\, (q^{20},q^{28};q^{48})_\infty}{(q;q)_\infty}\\
& = \frac{1}{(q;q)_\infty} \sum_{n=-\infty}^\infty q^{12n(3n+1)}\,(q^{-6n}-q^{6n+2}) 
\tag{By \eqref{eq:8}, with $q$ replaced by $q^{24}$ and $z$ replaced by $q^2$}\\
& = \frac{1}{(q;q)_\infty} \sum_{n=-\infty}^\infty (q^{6n(6n+1)}-q^{(6n+1)(6n+2)}) 
\end{align*}	
and
\begin{align*}
\sum_{n=0}^\infty d(n)\, q^n
& = \frac{(q^4,q^{10},q^{14},q^{24},q^{34},q^{38},q^{44},q^{48};q^{48})_\infty}{(q;q)_\infty}\\
& = \frac{(q^{10},q^{14},q^{24};q^{24})_\infty\, (q^{4},q^{44};q^{48})_\infty}{(q;q)_\infty}\\
& = \frac{1}{(q;q)_\infty} \sum_{n=-\infty}^\infty q^{12n(3n+1)}\,(q^{-30n}-q^{30n+10}) 
\tag{By \eqref{eq:8}, with $q$ replaced by $q^{24}$ and $z$ replaced by $q^{10}$}\\
& = \frac{1}{(q;q)_\infty} \sum_{n=-\infty}^\infty (q^{6n(6n-3)}-q^{(6n+2)(6n+5)}).
\end{align*}
We can write
\begin{align*}
& \sum_{n=0}^\infty \big(c(n)-d(n-2)\big)\, q^n \\
& = \sum_{n=0}^\infty c(n)\, q^n - \sum_{n=0}^\infty d(n)\, q^{n+2}\\
& = \frac{1}{(q;q)_\infty}  \sum_{n=-\infty}^\infty (q^{6n(6n+1)}-q^{(6n+1)(6n+2)} 
- q^{6n(6n-3)+2}+q^{(6n+2)(6n+5)+2})  \\
& = \frac{1}{(q;q)_\infty}  \sum_{n=-\infty}^\infty (q^{6n(6n+1)}-q^{(6n+1)(6n+2)} 
- q^{(6n-1)(6n-2)}+q^{(6n+3)(6n+4)}) \\
& = \frac{1}{(q;q)_\infty}  \sum_{n=-\infty}^\infty (q^{3n(3n-1)}-q^{(3n+1)(3n+2)}) \\
& = \frac{1}{(q;q)_\infty}  \sum_{n=-\infty}^\infty q^{3n(3n+1)}(q^{-6n}-q^{6n+2})  \\
&  = \frac{(q^2,q^4,q^6;q^6)_\infty\, (q^2,q^{10};q^{12})_\infty}{(q;q)_\infty} 
\tag{By \eqref{eq:8}, with $q$ replaced by $q^{6}$ and $z$ replaced by $q^{2}$}\\
&  = \frac{(q^2;q^2)_\infty\, (q^2,q^{10};q^{12})_\infty}{(q;q^2)_\infty\, (q^2;q^2)_\infty}\\
& =  (-q;q)_\infty\, (q^2,q^{10};q^{12})_\infty\\
& =  (-q;q^2)_\infty\, (-q^2;q^2)_\infty\, (q^2,q^{10};q^{12})_\infty\\
& =  (-q;q^2)_\infty\, (-q^6;q^6)_\infty\, (-q^2,-q^4;q^6)_\infty\, (q^2,q^{10};q^{12})_\infty\\
& =  (-q;q^2)_\infty\, (-q^6;q^6)_\infty\, \frac{(q^2,q^4,q^8,q^{10};q^{12})_\infty}{(q^2,q^4;q^6)_\infty}\\
& = \sum_{n=0}^\infty Q_2(n)\, q^n.
\end{align*}
and
\begin{align*}
& \sum_{n=0}^\infty \big(c(n)+d(n-2)\big)\, q^n \\
& = \frac{1}{(q;q)_\infty}  \sum_{n=-\infty}^\infty (q^{6n(6n+1)}-q^{(6n+1)(6n+2)} 
+ q^{(6n-1)(6n-2)}-q^{(6n+3)(6n+4)}) \\
& = \frac{1}{(q;q)_\infty}  \sum_{n=-\infty}^\infty (q^{6n(6n+1)}-q^{(6n+3)(6n+4)}) \\
& = \frac{1}{(q;q)_\infty}  \sum_{n=-\infty}^\infty (-1)^n q^{3n(3n-1)} \\
& = \frac{(q^6;q^6)_\infty}{(q;q)_\infty} \tag{By \eqref{eq:3} with $q$ replaced by $q^6$}\\
& = \sum_{n=0}^\infty b_6(n)\, q^n.
\end{align*}
This concludes the proof.

\section{Proof of Theorems \ref{Th:5} and \ref{Th:5a}}
\label{S4}

	\allowdisplaybreaks{
	G. E. Andrews and M. Merca \cite{Andrews12} proved the following truncated form of \eqref{eq:3}: For any $k\geqslant 1$, 
	\begin{equation} \label{eq:9}
	\frac{1}{(q;q)_\infty} \sum_{n=-(k-1)}^{k} (-1)^{n}\, q^{n(3n-1)/2}
	= 1+ (-1)^{k-1} \sum_{n=k}^\infty \frac{q^{{k\choose 2}+(k+1)n}}{(q;q)_n}
	\begin{bmatrix}
	n-1\\k-1
	\end{bmatrix},
	\end{equation}
	where 
	$$
	\begin{bmatrix}
	n\\k
	\end{bmatrix} 
	=
	\begin{cases}
	\dfrac{(q;q)_n}{(q;q)_k(q;q)_{n-k}}, &  \text{if $0\leqslant k\leqslant n$},\\
	0, &\text{otherwise.}
	\end{cases}
	$$
	We note that the series on the right hand side  of \eqref{eq:9} is the generating function for $M_k(n)$, i.e.,
	\begin{align*}
	\sum_{n=0}^\infty M_k(n)\, q^n = \sum_{n=k}^\infty \frac{q^{{k\choose 2}+(k+1)n}}{(q;q)_n}
	\begin{bmatrix}
	n-1\\k-1
	\end{bmatrix}. 
	\end{align*}
	By \eqref{eq:9}, with $q$ replaced by $q^6$, we get
	\begin{align*}
	\frac{1}{(q^6;q^6)_\infty} \sum_{n=-(k-1)}^{k} (-1)^{n}\, q^{3n(3n-1)}
	= 1+ (-1)^{k-1} \sum_{n=0}^\infty M_k(n)\, q^{6n}.
	\end{align*}
	Multiplying both sides of this identity by $$\frac{(q^6;q^6)_\infty}{(q;q)_\infty},$$ we obtain
	\begin{align*}
		\frac{1}{(q;q)_\infty} \sum_{n=-(k-1)}^{k} (-1)^{n}\, q^{3n(3n-1)}
	- \frac{(q^6;q^6)_\infty}{(q;q)_\infty}= (-1)^{k-1} \frac{(q^6;q^6)_\infty}{(q;q)_\infty} \sum_{n=0}^\infty M_k(n)\, q^{6n}
	\end{align*}
	or
	\begin{align*}
	& \left( \sum_{n=0}^\infty p(n)\, q^n \right) \left( \sum_{n=-(k-1)}^{k} (-1)^{n}\, q^{3n(3n-1)} \right)  - \sum_{n=0}^\infty b_6(n)\, q^n  \\
	& \qquad = (-)^{k-1} \left( \sum_{n=0}^\infty b_6(n)\, q^n \right) \left( \sum_{n=0}^\infty M_k(n)\, q^{6n} \right).
	\end{align*}
	The assertion of Theorem \ref{Th:5} follows by comparing coefficients of $q^n$ on both
	sides of this equation.

The proof of Theorem \ref{Th:5a} is quite similar to the proof of Theorem \ref{Th:5}.
In \cite{Xia}, E.~X.~W.~Xia and X.~Zhao considered Euler's pentagonal number theorem 
\eqref{eq:3}	
and they proved the following truncated form: For any $k\geqslant 1$, 
\begin{equation} \label{TPNTa}
\frac{1}{(q;q)_\infty} \sum_{n=-k}^{k} (-1)^{n}\, q^{n(3n-1)/2}= 1+(-1)^{k}\, \frac{q^{k(k+1)/2}}{(q;q)_k} \sum_{n=0}^\infty \frac{q^{(n+k+1)(k+1)}}{(q^{n+k+1};q)_\infty}.
\end{equation}
We remark that the series on the right hand side  of \eqref{TPNTa} is the generating function for $\widetilde{P}_k(n)$, i.e.,
\begin{equation*}
\sum_{n=0}^\infty \widetilde{P}_k(n)\, q^n = \frac{q^{k(k+1)/2}}{(q;q)_k} \sum_{n=0}^\infty \frac{q^{(n+k+1)(k+1)}}{(q^{n+k+1};q)_\infty}. 
\end{equation*}
	By \eqref{TPNTa}, with $q$ replaced by $q^6$, we get
\begin{align*}
\frac{1}{(q^6;q^6)_\infty} \sum_{n=-k}^{k} (-1)^{n}\, q^{3n(3n-1)}= 1+(-1)^{k}\, \frac{q^{3k(k+1)}}{(q^6;q^6)_k} \sum_{n=0}^\infty \frac{q^{6(n+k+1)(k+1)}}{(q^{6(n+k+1)};q^6)_\infty}.
\end{align*}
Multiplying both sides of this identity by the generating function of $b_6(n)$,
we obtain
\begin{align*}
(-1)^{k}\left( \Big(\sum_{n=1}^\infty p(n)\, q^n \Big) \Big( \sum_{n=-k}^k (-1)^n\, q^{n(3n-1)/2} 
\Big) - \sum_{n=1}^\infty b_6(n)\, q^{n}  \right)  \\
= \left( \sum_{n=1}^\infty b_6(n)\, q^{n} \right) \left( \sum_{n=0}^\infty \widetilde{P}_k(n)\, q^{6n}\right).
\end{align*}
The proof of Theorem \ref{Th:5a} follows easily considering Cauchy's multiplication of two power series.

\section{Proof of Theorem \ref{Th:6}}
\label{S5}

\allowdisplaybreaks{
The Jacobi triple product identity (cf.\ \cite[Eq.~(1.6.1)]{GaRaAA}) states that
\begin{equation} \label{eq:10} 
\sum_{n=-\infty}^\infty (-z)^n\, q^{n(n-1)/2}=(z,q/z,q;q)_\infty.
\end{equation} 
Considering \eqref{eq:10} with $q$ replaced by $q^6$ and $z$ replaced by $q$, we can write
\begin{align*}
\sum_{n=-\infty}^\infty (-1)^n\, q^{n(3n-2)} 
& = (q,q^5,q^6;q^6)_\infty \\
& = \frac{(q^6;q^6)_\infty}{(q;q)_\infty} \cdot (q,q^5;q^6)_\infty\, (q,q^2,q^3;q^3)_\infty\\
& = \frac{(q^6;q^6)_\infty}{(q;q)_\infty} \sum_{n=-\infty}^\infty q^{3n(3n+1)/2}\, (q^{-3n}-q^{3n+1})
\tag{By \eqref{eq:8}, with $q$ replaced by $q^3$ and $z$ replaced by $q$}\\
& = \frac{(q^6;q^6)_\infty}{(q;q)_\infty} \left(  \sum_{n=-\infty}^\infty q^{3n(3n-1)/2} - \sum_{n=-\infty}^\infty q^{(3n+1)(3n+2)/2} \right) \\
& = \frac{(q^6;q^6)_\infty}{(q;q)_\infty} \left(  \sum_{\substack{n=0\\n\not\equiv 1 \pmod 3}}^\infty q^{n(n+1)/2} - \sum_{\substack{n=0\\n\equiv 1 \pmod 3}}^\infty 2\,q^{n(n+1)/2} \right) \\
& = \left( \sum_{n=0}^\infty b_6(n)\, q^n  \right) \left( \sum_{n=0}^\infty \rho_n\, q^{n(n+1)/2}\right)\\
& = \sum_{n=0}^\infty \left( \sum_{j=0}^n \rho_j\, b_6\big(n-j(j+1)/2\big)\right) q^n.
\end{align*}
This concludes the proof.
}

\section{Proof of Theorems \ref{Th:8} and \ref{Th:10}}
\label{S6}

\subsection{Analytic proof}

\allowdisplaybreaks{
We have
\begin{align}
\sum_{n=0}^\infty Q_2(n)\, q^n 
& = (-q;q^2)_\infty\, (-q^6;q^6)_\infty \nonumber \\
& = \frac{(q^2;q^4)_\infty}{(q;q^2)_\infty}\, \frac{1}{(q^6;q^{12})_\infty}\nonumber \\
& = \frac{1}{(q;q)_\infty} \cdot (q^2,q^{10};q^{12})_\infty\, (q^2,q^4,q^6;q^6)_\infty\nonumber \\
& = \frac{1}{(q;q)_\infty} \sum_{n=-\infty}^\infty q^{3n(3n+1)}\, (q^{-6n}-q^{6n+2})
\tag{By \eqref{eq:8}, with $q$ replaced by $q^6$ and $z$ replaced by $q^2$}\nonumber \\
& = \frac{1}{(q;q)_\infty} \left(  \sum_{n=-\infty}^\infty q^{3n(3n-1)} - \sum_{n=-\infty}^\infty q^{(3n+1)(3n+2)} \right) \nonumber \\
& = \frac{1}{(q;q)_\infty} \left(  \sum_{\substack{n=0\\n\not\equiv 1 \pmod 3}}^\infty q^{n(n+1)} - \sum_{\substack{n=0\\n\equiv 1 \pmod 3}}^\infty 2\,q^{n(n+1)} \right)\nonumber \\
& = \left( \sum_{n=0}^\infty p(n)\, q^n  \right) \left( \sum_{n=0}^\infty \rho_n\, q^{n(n+1)}\right)\nonumber \\
& = \sum_{n=0}^\infty \left( \sum_{j=0}^n \rho_j\, p\big(n-j(j+1)\big)\right) q^n\label{eq:11}
\end{align}
and
\begin{align*}
\sum_{n=0}^\infty Q_2(n)\, q^n
& = (-q^3;q^3)_\infty\, (-q,-q^5;q^6)_\infty\\
& = \frac{(-q^3;q^3)_\infty}{(q^6;q^6)_\infty} \cdot (-q,-q^5,q^6;q^6)_\infty\\
& = \frac{1}{(q^3;q^3)_\infty} \sum_{n=-\infty}^\infty q^{n(3n-2)} 
\tag{By  \eqref{eq:10}, with $q$ replaced by $q^6$ and $z$ replaced by $-q$}\\
& = \left( \sum_{n=0}^\infty p(n)\, q^{3n} \right) \left( \sum_{n=-\infty}^\infty q^{n(3n-2)} \right) \\
& = \sum_{n=0}^\infty \left( \sum_{j=-\infty}^\infty p\left(\frac{n-j(3j-2)}{3} \right)  \right) q^n,
\end{align*} 
where $p(x)=0$ if $x$ is not a nonnegative integer. Theorem \ref{Th:8} is proved.

On the other hand, by these relations, we deduce that
\begin{align*}
(q;q)_\infty \sum_{n=0}^\infty Q_2(n)\, q^n  = \sum_{n=0}^\infty \rho_n\, q^{n(n+1)}
\end{align*}
and
\begin{align*}
(q^3;q^3)_\infty \sum_{n=0}^\infty Q_2(n)\, q^n = \sum_{n=-\infty}^\infty q^{n(3n-2)}.
\end{align*}
Considering \eqref{eq:3}, these equations can be rewritten as
\begin{align*}
\left( \sum_{n=-\infty}^\infty (-1)^n\, q^{n(3n-1)/2}\right) \left(  \sum_{n=0}^\infty Q_2(n)\, q^n \right)  = \sum_{n=0}^\infty \rho_n\, q^{n(n+1)}
\end{align*}
and
\begin{align*}
\left( \sum_{n=-\infty}^\infty (-1)^n\, q^{3n(3n-1)/2}\right) \left( \sum_{n=0}^\infty Q_2(n)\, q^n \right)  = \sum_{n=-\infty}^\infty q^{n(3n-2)}.
\end{align*}
The assertions of Theorem \ref{Th:10} follows easily by comparing coefficients of $q^n$ on both
sides of these equations.
}

\subsection{Combinatorial proof of Theorem \ref{Th:8} (ii)}


We will use a particular case of \cite[Lemma 2.1]{BM22}. We first introduce some notation. 
 Let $\mathcal Q_{6,1}(n)$ be the set of partitions of $n$ into distinct parts congruent to $\pm 1\pmod{6})$. Let $\mathcal W_{6,1}(n)$ be the set of pairs $(\mu, k(3k-2))$, where $\mu$ is a partition into parts divisible by $6$,  $k\in \mathbb Z$, and $|\mu|+ k(3k-2)=n$. Then, \cite[Lemma 2.1]{BM22} with $m=6$ and $r=1$ gives a bijection $\xi_{6,1}: \mathcal Q_{6,1}(n)\to \mathcal W_{6,1}(n)$.

Let $n$ be a nonnegative integer. We create a  bijection $$\psi: \mathcal Q_2(n)\to \bigcup_{k\in \mathbb Z} \mathcal P\left(\frac{n-j(3j-2)}{3}\right).$$ 
Start with $\lambda\in \mathcal Q_2(n)$. Write $\lambda=\alpha\cup \beta$,  where $\alpha$ is a  partition into distinct parts congruent to $1$ or $5$ modulo $6$,  and $\beta$ is a partition into distinct parts divisible by $3$. Thus, $\beta_{/3}$ is a  partition with distinct parts and $\varphi_{Gl}^{-1}(\beta_{/3})$ is a partition  with odd parts. Moreover, $3\varphi_{Gl}^{-1}(\beta_{/3})$ is a partition whose  parts are  congruent to $3$ modulo $6$ and $|3\varphi_{Gl}^{-1}(\beta_{/3})|=|\beta|$. 
Let $\xi_{6,1}(\alpha)=(\mu, k(3k-2))$ for some $k\in \mathbb Z$. Since the parts of $\mu$ are divisible by $6$, all parts of $\mu\cup 3\varphi_{Gl}^{-1}(\beta_{/3})$ are divisible by $3$. Define $\psi(\lambda)=(\mu\cup3\varphi_{Gl}^{-1}(\beta_{/3}))_{/3}$. Since $|\mu|=|\alpha|-k(3k-2)$ and $|\alpha|+|\beta|=n$, it follows that $\psi(\lambda)$ is a partition of $\frac{n-k(3k-2)}{3}$.

For the inverse, let $k \in \mathbb Z$. Start with a partition $\eta$ of $\frac{n-k(3k-2)}{3}$. Then $3\eta$ is a partition of $ n-k(3k-2)$. Write $3\eta=\mu\cup\pi$, where $\mu$ (respectively $\pi$) has parts congruent to $0$ (respectively $3$) modulo $6$. We have that  $\xi_{6,1}^{-1}(\mu, k(3k-2))$ is a partition of $|\mu|+k(3k-2)$ into parts congruent to $1$ or $5$ modulo $6$. The partition $\pi/3$ has odd parts and  $3\varphi_{Gl}(\pi_{/3})$ is a partition into distinct parts divisible by  $3$ and $|3\varphi_{Gl}(\pi/3)|=|\pi|$. Then $\psi^{-1}(\eta)=\xi_{6,1}^{-1}(\mu, k(3k-2)) \cup 3\varphi_{Gl}(\pi/3)\in \mathcal Q_2(n)$.


\subsection{Combinatorial proof of Theorem \ref{Th:10} (ii)}

We denote by $\mathcal Q(n)$ the set of  partitions of $n$ into distinct parts and we set $\mathcal Q:=\cup_{n\geq 0} \mathcal Q(n)$,  $\mathcal Q_2:=\cup_{n\geq 0} \mathcal Q_2(n)$. and $\mathcal Q_{6,1}:=\cup_{n\geq 0} \mathcal Q_{6,1}(n)$.  Let $\varphi_F$ be the involution defined by Franklin to give a combinatorial proof of Euler's pentagonal number theorem (see, for example, \cite[Theorem 1.6]{Andrews98}). Let $$\mathcal A(n):=\{ (\lambda, \mu): \lambda\in \mathcal Q_2, \mu =3\eta  \text{ with } \eta\in \mathcal Q, |\lambda|+|\mu|=n\}$$ and $$\mathcal{EA}(n):=\{(\lambda, \mu)\in \mathcal A(n): \mu_{/3} \text{ pentagonal partition}\}.$$ Here a pentagonal partition is either $\emptyset$ or a partition of the form $(2i, 2i-1, \ldots, i+1)$ or $(2i-1, 2i-2, \ldots, i)$ for some integer $i>0$.

The involution $(\lambda, \mu)\mapsto (\lambda, 3\varphi_F(\mu_{/3}))$ on $\mathcal A(n)\setminus \mathcal{EA}(n) $ proves combinatorially that $$\displaystyle{\sum_{j=-\infty}^\infty (-1)^j\, Q_2\big(n-3j(3j-1)/2 \big)}$$ is the generating function for \begin{align}\label{A-diff} |\{(\lambda, \mu) \in \mathcal A(n): \ell(\mu) \text{ even}\}|-|\{(\lambda, \mu) \in \mathcal A(n): \ell(\mu) \text{ odd}\}|. \end{align}

We define another involution on a subset of $\mathcal A(n)$ that reverses the parity of the length of the second partition in the pair. Given  $(\lambda, \mu)\in \mathcal A(n)$, we write $\lambda=\lambda^{3\mid}\cup \lambda^{3\nmid}$, where the parts of $\lambda^{3\mid}$ are all parts of $\lambda$ which are divisible by $3$. If $\ell(\lambda^{3\mid})\not \equiv \ell(\mu)\pmod 2$, we map $(\lambda^{3\mid}, \lambda^{3\nmid}, \mu)$ to $(\mu, \lambda^{3\nmid}, \lambda^{3\mid})$. If $\ell(\lambda^{3\mid})\not \equiv \ell(\mu)\pmod 2$ and $\lambda^{3\mid}\neq \mu$, let $i$ be the smallest integer such that $\lambda^{3\mid}_i\neq \mu_i$. If $\lambda^{3\mid}_i>\mu_i$, we remove part $\lambda^{3\mid}_i$ from $\lambda^{3\mid}$ and insert a part equal to $\lambda^{3\mid}_i$ into $\mu$. If $\lambda^{3\mid}_i<\mu_i$, we remove part $\mu_i$ from $\mu$ and insert a part equal to $\mu_i$ into $\lambda^{3\mid}$. We obtain an involution on the set $\{(\lambda, \mu)\in \mathcal A(n): \lambda^{3\mid} \neq \mu\}$ that reverses the parity of $\ell(\mu)$. Thus, \eqref{A-diff} equals  \begin{align*} |\{(\lambda, \lambda^{3\mid}) \in \mathcal A(n): \ell(\lambda^{3\mid}) \text{ even}\}|-|\{(\lambda, \lambda^{3\mid}) \in \mathcal A(n): \ell(\lambda^{3\mid}) \text{ odd}\}|. \end{align*} Mapping $(\lambda, \lambda^{3\mid})=(\lambda^{3\mid}, \lambda^{3\nmid}, \lambda^{3\mid})$ to $(\lambda^{3\nmid}, 2\lambda^{3\mid})$ and setting $$\mathcal B(n):=\left\{(\alpha, \beta): \alpha\in \mathcal Q_{6,1},  \beta =6\gamma  \text{ with } \gamma\in \mathcal Q,  |\alpha|+|\beta|=n\right\},$$ we see that \eqref{A-diff} equals \begin{align*} |\{(\alpha, \beta)\in \mathcal B(n): \ell(\beta) \text{ even}\}|-|\{(\alpha, \beta)\in \mathcal B(n): \ell(\beta) \text{ odd}\}|. \end{align*}

Let $\mathcal C(n)$ be the set of triples $(\gamma, k(3k-2), \beta)$ such that $k\in \mathbb Z$, $\gamma$ and $\beta$ are partitions with parts divisible by $6$, $\beta\in \mathcal Q$, and $|\gamma|+k(3k-2)+|\beta|=n$.

We define a bijection from $\mathcal B(n)$ to the set $\mathcal C(n)$  by $$(\alpha, \beta)\mapsto (\xi_{6,1}(\alpha), \beta)=(\gamma, k(3k-2), \beta)$$ where, $\xi_{6,1}: \mathcal Q_{6,1}(n)\to \mathcal W_{6,1}(n)$ is the bijection of \cite[Lemma 2.1]{BM22}.
Then, \eqref{A-diff} equals \begin{align*} |\{(\gamma, k(3k-2), \beta)\in \mathcal C(n): \ell(\beta) \text{ even}\}|-|\{(\gamma, k(3k-2), \beta)\in \mathcal C(n): \ell(\beta) \text{ odd}\}|. \end{align*}

Finally, we define an involution $\zeta$ on $$\{(\gamma, k(3k-2), \beta)\in C(n):(\gamma, \beta)\neq (\emptyset, \emptyset)\}.$$  Start with $(\gamma, k(3k-2), \beta)\in \mathcal C(n)$ with $(\gamma, \beta)\neq (\emptyset, \emptyset)$. If $\gamma_1>\beta_1$, remove part $\gamma_1$ from $\gamma$ and insert a part equal to $\gamma_1$ into $\beta$. If $\gamma_1\leq \beta_1$, remove part $\beta_1$ from $\beta$ and insert a part equal to $\beta_1$ into $\gamma$. The involution $\zeta$ reverses the parity of $\ell(\beta)$.

This completes the combinatorial proof of Theorem \ref{Th:10}(ii).

\section{Inequalities and open problems}
\label{S7}

Linear inequalities involving partition functions, especially Euler's partition function $p(n)$, have been the subject of recent studies 
by
G. E. Andrews and M. Merca \cite{Andrews12,Andrews18,Andrews20},
C. Ballantine and M. Merca \cite{BM,BMaa,BMa},
C. Ballantine, M. Merca, D. Passary and A. J. Yee \cite{BMPY},
V. J. W. Guo and J. Zeng \cite{Guo},
J. Katriel \cite{Katriel},
M. Merca \cite{Merca12,Merca15,Merca16,Merca18,Merca19c,Merca19bb,Merca19b,Merca21b,Merca21c,Merca21aa,Merca21a,Merca21d,Merca21e,Merca21f},
M. Merca and J. Katriel \cite{MK},
M. Merca, C. Wang and A. J. Yee \cite{MWY},
M. Merca and A. J. Yee \cite{MY}.
For example,  G. E. Andrews and M. Merca \cite{Andrews12} proved that: for $n\geqslant 0$, $k>0$,
\begin{align*}
(-1)^{k-1} \sum_{j=-(k-1)}^{k} (-1)^{j}\, p\big(n-j(3j-1)/2\big) \geqslant 0.
\end{align*} 
Recently \cite[Corollary 11]{Andrews18}, the same authors found a new infinite family of linear homogeneous inequalities for $p(n)$ which involves the triangular numbers: if at least one of $n$ and $k$ is odd,
\begin{align*}
(-1)^{k-1} \sum_{j=0}^{2k-1} (-1)^{j(j-1)/2}\, p\big(n-j(j+1)/2\big) \geqslant 0.
\end{align*}
As a consequence of Theorem \ref{Th:5}, we remark a new infinite family of linear inequalities for  $p(n)$.

\begin{corollary}
	For $n\geqslant 0$, $k>0$,
	\begin{align*}
(-1)^{k} \left(b_6(n) - \sum_{j=-(k-1)}^k (-1)^j\, p\big(n-3j(3j-1)\big) \right)
\geqslant 0, 
\end{align*}
with strict inequality if $n\geqslant 3k(3k+1)$. 
\end{corollary}

For example, the cases $k=1$ and $k=2$ of this corollary provides the following double inequality
\begin{align}\label{C8_1_2}
p(n)-p(n-6)-p(n-12)+p(n-30) \leqslant b_6(n) \leqslant p(n)-p(n-6).
\end{align}
In terms of $\mex_{2,2}$, inequality \eqref{C8_1_2} becomes $$pm_2(n)+pm_4(n)-pm_8(n)-pm_{10}(n) \leqslant b_6(n) \leqslant pm_2(n)+pm_4(n).$$

%
%

In this section, inspired by the identity \eqref{eq:5} and Theorems \ref{Th:6}, \ref{Th:8} and \ref{Th:10}, we propose as conjectures  new infinite families of linear inequalities for the partition functions $p(n)$, $b_6(n)$ and $Q_2(n)$.

\subsection{Euler's partition function}

Inspired by Theorem \ref{Th:8}, for $k\geqslant 0$ we investigated the following series:
 \begin{align*}
 (-q;q^2)_\infty\, (-q^6;q^6)_\infty - \frac{1}{(q;q)_\infty} \sum_{j=0}^k \rho_j\, q^{j(j+1)}.
 \end{align*} 
There is a substantial amount of numerical evidence to conjecture 
that this series  has nonnegative coefficients if $k$ is not congruent to $0$ modulo $3$ and nonpositive coefficients if $k$ is congruent to $0$ modulo $3$.
In addition, we conjecture that the coefficient of $q^n$ in this series is nonzero if and only if $n\geqslant (k+1)(k+2)$. We have the following combinatorial interpretation of this conjecture.

\begin{conjecture}
	For $n,k\geqslant 0$,
	\begin{enumerate}
		\item[(i)] $\displaystyle{\sum_{j=0}^{3k}\rho_j\, p\big(n-j(j+1)\big) \geqslant Q_2(n)} $,
		\item[] with strict inequality if $n\geqslant (3k+1)(3k+2)$;
		\item[(ii)] $\displaystyle{\sum_{j=0}^{3k+1}\rho_j\, p\big(n-j(j+1)\big) \leqslant Q_2(n)} $;
		\item[] with strict inequality if $n\geqslant (3k+2)(3k+3)$.
	\end{enumerate}
\end{conjecture}

Assuming this conjecture, we remark the following double inequality:
\begin{align}
& p(n)-2p(n-2)+p(n-6) \leqslant Q_2(n) \notag \\
& \qquad\qquad\qquad \leqslant p(n)-2p(n-2)+p(n-6)+p(n-12). \label{conj2_1}
\end{align}

In terms of $\mex_{2,2}$, inequality \eqref{conj2_1} becomes $$pm_2(n)-pm_4(n) \leqslant Q_2(n) \leqslant pm_2(n)-pm_4(n)+pm_{>6}(n).$$

\subsection{$6$-regular partitions}

Inspired by the identity \eqref{eq:5} and the truncated pentagonal number theorem \eqref{eq:9}, for $k>0$ we considered the following series:
\begin{align*}
(q^6;q^6)_\infty - \frac{(q^6;q^6)_\infty}{(q;q)_\infty} \sum_{n=-(k-1)}^{k} (-1)^{n}\, q^{n(3n-1)/2}.
\end{align*}
There is a substantial amount of numerical evidence to conjecture 
that this series has nonnegative coefficients if $k$ is even and nonpositive coefficients if $k$ is odd.
In addition, we conjecture that the coefficient of $q^n$ in this series is nonzero if and only if $n\geqslant k(3k+1)/2$. We have the following combinatorial interpretation of this conjecture. For any integer $n$, let
$$
\alpha_n := \begin{cases}
(-1)^m,& \text{if $n=3m(3m-1)$, $m\in\mathbb{Z},$}\\
0, & \text{otherwise.}
\end{cases}
$$

\begin{conjecture}\label{CJ2}
	For $n\geqslant 0$, $k>0$,
	\begin{align*}
	(-1)^{k} \left( \alpha_n - \sum_{j=-(k-1)}^k (-1)^j\, b_6\big(n-j(3j-1)/2\big)  \right)  \geqslant 0,
	\end{align*}
	with strict inequality if $n\geqslant k(3k+1)/2$
\end{conjecture}

We remark that this inequality can be rewritten in terms of $M_k(n)$ as
follows: for $n\geqslant 0$, $k>0$,
\begin{align}
\sum_{j=-\infty}^\infty (-1)^j\, M_k\big(n-3j(3j-1) \big) \geqslant 0,	\label{eq:12}
\end{align}
with strict inequality if $n\geqslant k(3k+1)/2$.

In analogy with Conjecture \ref{CJ2}, 
we also make the following conjecture.

\begin{conjecture}\label{CJ2a}
	For $n\geqslant 0$, $k>0$,
	\begin{align*}
	(-1)^{k-1} \left( \alpha_n - \sum_{j=-k}^k (-1)^j\, b_6\big(n-j(3j-1)/2\big)  \right)  \geqslant 0,
	\end{align*}
	with strict inequality if $n\geqslant k(3k+1)/2$.
\end{conjecture}

It is easy to see that Conjecture \ref{CJ2a} is a weaker version of Conjecture \ref{CJ2}.
The inequality given by Conjecture \ref{CJ2a} can be rewritten in terms of $\widetilde{P}_k(n)$ as
follows: for $n\geqslant 0$, $k>0$,
\begin{align}
\sum_{j=-\infty}^\infty (-1)^j\, \widetilde{P}_k\big(n-3j(3j-1) \big) \geqslant 0,	\label{eq:12a}
\end{align}
with strict inequality if $n\geqslant k(3k+1)/2$. Clearly the inequality \eqref{eq:12} implies the inequality \eqref{eq:12a}. 

Inspired by Theorem \ref{Th:6}, for $k\geqslant 0$ we investigated the following series:
\begin{align*}
(q,q^5,q^6;q^6)_\infty - \frac{(q^6;q^6)_\infty}{(q;q)_\infty} \sum_{j=0}^k \rho_j\, q^{j(j+1)/2}.
\end{align*}
There is a substantial amount of numerical evidence to conjecture 
that this series has nonnegative coefficients if $k$ is not congruent to $0$ modulo $3$ and nonpositive coefficients if $k$ is congruent to $0$ modulo $3$.
In addition, the coefficient of $q^n$ in this series is nonzero if and only if $n\geqslant (k+1)(k+2)/2$. We have the following combinatorial interpretation of this conjecture. For any integer $n$, let 
$$
\beta_n := \begin{cases}
(-1)^m,& \text{if $n=m(3m-2)$, $m\in\mathbb{Z},$}\\
0, & \text{otherwise.}
\end{cases}
$$

\begin{conjecture}
	For $n,k\geqslant 0$,
	\begin{enumerate}
		\item[(i)] $\displaystyle{\sum_{j=0}^{3k}\rho_j\, b_6\big(n-j(j+1)/2\big) \geqslant \beta_n} $,
		\item[] with strict inequality if $n\geqslant (3k+1)(3k+2)/2$;
		\item[(ii)] $\displaystyle{\sum_{j=0}^{3k+2}\rho_j\, b_6\big(n-j(j+1)/2\big) \leqslant \beta_n} $;
		\item[] with strict inequality if $n\geqslant (3k+2)(3k+3)/2$.
	\end{enumerate}
\end{conjecture}

Assuming this conjecture, we remark the following double inequality:
\begin{align*}
b_6(n)-2b_6(n-1)+b_6(n-3) \leqslant \beta_n \leqslant b_6(n)-2b_6(n-1)+b_6(n-3)+b_6(n-6).
\end{align*}

\subsection{Partitions into distinct parts $\not\equiv \pm 2 \pmod 6$}

Inspired by Theorem \ref{Th:10}.(i), for $k>0$ we considered the following series:
\begin{align*}
	(q^2,q^{10};q^{12})_\infty\, (q^2;q^2)_\infty - (-q;q^2)_\infty\, (-q^6;q^6)_\infty  \sum_{j=-(k-1)}^k (-1)^j q^{j(3j-1)/2}.
\end{align*}

There is a substantial amount of numerical evidence to conjecture 
that this series has nonnegative coefficients if $k$ is even and nonpositive coefficients if $k$ is odd.
We have the following combinatorial interpretation of this conjecture. For any nonnegative integer $n$, let 
$$
\gamma_n := \begin{cases}
	\rho_m,& \text{if $n=m(m+1)$, $m\in\mathbb{N}_0,$}\\
	0, & \text{otherwise.}
\end{cases}
$$

\begin{conjecture}\label{CJ5}
	For $n\geqslant 0$, $k>0$,
	\begin{align*}
		(-1)^{k} \left(\gamma_n - \sum_{j=-(k-1)}^k (-1)^j\, Q_2\big(n-j(3j-1)/2\big)  \right)  \geqslant 0.
	\end{align*}
\end{conjecture}

We remark that this inequality can be rewritten in terms of $M_k(n)$ as
follows: for $n\geqslant 0$, $k>0$,
\begin{align}
\sum_{j=0}^\infty \rho_j\, M_k\big(n-j(j+1) \big) \geqslant 0.	\label{eq:13}
\end{align}

 We  also make the following conjecture which is 
  weaker than Conjecture \ref{CJ5}.

\begin{conjecture}\label{CJ5a}
	For $n\geqslant 0$, $k>0$,
	\begin{align*}
	(-1)^{k-1} \left(\gamma_n - \sum_{j=-k}^k (-1)^j\, Q_2\big(n-j(3j-1)/2\big)  \right)  \geqslant 0.
	\end{align*}
\end{conjecture}

This inequality can be rewritten in terms of $\widetilde{P}_k(n)$ as
follows: for $n\geqslant 0$, $k>0$,
\begin{align}
\sum_{j=0}^\infty \rho_j\, \widetilde{P}_k\big(n-j(j+1) \big) \geqslant 0.	\label{eq:13a}
\end{align}
It is clear that the inequality \eqref{eq:13} implies the inequality \eqref{eq:13a}.

Regarding inequalities \eqref{eq:12} - \eqref{eq:13a}, it would be very appealing to have  combinatorial interpretations for their sums.

\bigskip



\begin{thebibliography}{00}


\bibitem{Ahmed}
Z.~Ahmed, N.~D.~Baruah,
New congruences for $\ell$-regular partitions for $\ell\in \{5, 6, 7, 49\}$,
\textit{Ramanujan J} \textbf{40}(3) (2016) 649--668. 

\bibitem{Andrews98} 
G.~E.~Andrews, 
\textit{The Theory of Partitions}, 
Cambridge Mathematical Library, Cambridge University Press, Cambridge, 1998. Reprint of the 1976 original.

\bibitem{Andrews10}
G.~E.~Andrews,  M.~D.~Hirschhorn, J.~A.~Sellers, 
Arithmetic properties of partitions with even parts distinct,
\textit{Ramanujan J} \textbf{23}(1–3) (2010) 169--181. 

\bibitem{Andrews12} 
G.~E.~Andrews, M.~Merca,  The truncated pentagonal number theorem, 
\textit{J. Combin. Theory Ser. A}, \textbf{119} (2012) 1639--1643.

\bibitem{Andrews18} 
G.~E. Andrews, M.~Merca,  Truncated theta series and a problem of Guo and Zeng, 
\textit{J. Combin. Theory Ser. A}, \textbf{154} (2018) 610--619.

\bibitem{Andrews20} 
G.~E.~Andrews, M.~Merca,  
On the number of even parts in all partitions of $n$ into distinct parts, 
\textit{Ann. Comb.} \textbf{24} (2020) 47--54.

\bibitem{AN} G.~E.~Andrews, D.~Newman,  The minimal excludant in integer partitions. \textit{J. Integer Seq.} \textbf{23} (2020), no. 2, Art. 20.2.3, 11 pp.

\bibitem{BM}
C.~Ballantine, M.~Merca,
On identities of Watson type,
\textit{Ars Math. Contemp.} \textbf{17} (2019) 277--290.

\bibitem{BMaa}
C.~Ballantine, M.~Merca, 
The minimal excludant and colored partitions,
\textit{S\'{e}m. Lothar. Combin.} \textbf{84B} (2020) Article \#23.

\bibitem{BMa}
C.~Ballantine, M.~Merca,
Combinatorial proof of the minimal excludant theorem,
\textit{Int. J. Number Theory} \textbf{17}(8) (2021) 1765--1779.

\bibitem{BM22}
 C.~Ballantine, M.~Merca,  Almost 3-regular overpartitions. \textit{Ramanujan J.} \textbf{58} (2022), no. 3, 957--971. 

\bibitem{BMPY}
C.~Ballantine, M.~Merca, D.~Passary, A.~J.~Yee, 
Combinatorial proofs of two truncated theta series theorems,
\textit{J. Combin. Theory Ser. A} \textbf{160} (2018) 168--185.

\bibitem{BMR}C.~Ballantine, M.~Merca, C.-S.~Radu, Parity of $3$-regular partition numbers and Diophantine equations,  	arXiv:2212.09810, 2022

\bibitem{BW}
C.~Ballantine, A.~Welch, Amanda, PED and POD partitions: Combinatorial proofs of recurrence relations, \textit{Discrete Math.} \textbf{346} (2023), no. 3, Paper No. 113259.

\bibitem{Carlitz} 
L.~Carlitz, M.~V.~Subbarao, 
A simple proof of the quintuple product identity, 
\textit{Proc. Amer. Math. Soc.} \textbf{32} (1972), 42--44.

\bibitem{Carlson}
R.~Carlson, J.~J.~Webb,  
Infinite families of infinite families of congruences for $k$-regular partitions, 
\textit{Ramanujan J} \textbf{33} (2014) 329--337.

\bibitem{Cui}
S.-P.~Cui, N.~S.~S.~Gu,
Arithmetic properties of $\ell$-regular partitions,
\textit{Adv. Appl. Math.} \textbf{51}(4) (2013) 507--523.

\bibitem{Dand}
B.~Dandurand, D.~Penninston,
$\ell$-Divisibility of $\ell$-regular partition functions,
\textit{Ramanujan J} \textbf{19} (2009) 63--70.


\bibitem{Furcy}
D.~Furcy, D.~Penniston,
Congruences for $\ell$-regular partition functions modulo $3$,
\textit{Ramanujan J} \textbf{27} (2012) 101--108.
 
\bibitem{GaRaAA}
G.~Gasper, M.~Rahman, {\em Basic Hypergeometric Series},
Encyclopedia of Mathematics And Its Applications 35, Cambridge
University Press, Cambridge, 1990.

\bibitem{Guo} 
V.~J.~W.~Guo, J.~Zeng, Two truncated identities of Gauss, 
\textit{J. Combin. Theory Ser. A}, \textbf{120} (2013) 700--707.

\bibitem{Hirschhorn}
 M.~D.~Hirschhorn, J.~A.~Sellers, 
 Elementary proofs of parity results for $5$-regular partitions,
 \textit{Bull. Aust. Math. Soc.} \textbf{81}(1) (2010) 58--63.

\bibitem{Hou}
Q.-H.~Hou, L.~H.~Sun,  L.~Zhang,
Quadratic forms and congruences for $\ell$-regular partitions modulo $3$, $5$ and $7$,
\textit{Adv. Appl. Math.} \textbf{70} (2015) 32--44.

\bibitem{Katriel}
J.~Katriel,
Asymptotically trivial linear homogeneous partition inequalities,
\textit{J. Number Theory} \textbf{184} (2018) 107--121


\bibitem{Lovejoy}
J.~Lovejoy, D.~Penniston,
$3$-regular partitions and a modular $K3$ surface, 
\textit{Contemp. Math.} \textbf{291} (2001) 177--182.

\bibitem{Merca12} 
M.~Merca, 
Fast algorithm for generating ascending compositions, 
\textit{J. Math. Model. Algorithms} \textbf{11} (2012) 89--104. 

\bibitem{Merca15}
M.~Merca,
A generalization of Euler's pentagonal number recurrence for the partition function,
\textit{Ramanujan J} \textbf{37} (2015) 589--595.

\bibitem{Merca16} 
M.~Merca, 
A new look on the truncated pentagonal number theorem,
\textit{Carpathian J. Math.} \textbf{32} (2016) 97--101. 

\bibitem{Merca18} 
M.~Merca, 
Higher-order difference and higher-order partial sums of Euler's partition function,
\textit{Ann. Acad. Rom. Sci. Ser. Math. Appl.} \textbf{10}(1) (2018) 59--71. 

\bibitem{Merca19c}
M.~Merca,
On a combinatorial interpretation of the bisectional pentagonal number theorem
\textit{J. of Ramanujan Society of Mathematics and Mathematical Sciences} \textbf{7}(1) (2019) 7--18.

\bibitem{Merca19bb}
M.~Merca,
Truncated theta series and Rogers-Ramanujan functions,
\textit{Exp. Math.} \textbf{30}(3) (2021) 364--371.

\bibitem{Merca21b}
M.~Merca,
On the sum of parts in the partitions of $n$ into distinct parts,
\textit{Bull. Aust. Math. Soc.} \textbf{104}(2) (2021) 228--237.

\bibitem{Merca21c}
M.~Merca,
On the partitions into distinct parts and odd parts,
\textit{Quaest. Math.} \textbf{44}(8) (2021) 1095--1105.	

\bibitem{Merca21aa}
M.~Merca,
Generalized Lambert series and Euler's pentagonal number theorem,
\textit{Mediterr. J. Math.} \textbf{18}:29 (2021).

\bibitem{Merca21a}
M.~Merca,
Polygonal numbers and Rogers-Ramanujan-Gordon theorem,
\textit{Ramanujan J} \textbf{55}  (2021) 783--792.

\bibitem{Merca21d}
M.~Merca,
On the number of partitions into parts not congruent to $0, \pm 3 \pmod {12}$,
\textit{Period. Math. Hungar.} \textbf{83} (2021) 133--143.

\bibitem{Merca21e}
M.~Merca,
Linear inequalities concerning partitions into distinct parts,
\textit{Ramanujan J} \textbf{58}  (2022) 491--503.	

\bibitem{Merca19b}
M.~Merca,
On two truncated quintuple series theorems,
\textit{Exp. Math.} \textbf{31}(2) (2022) 606--610.

\bibitem{Merca21f}
M.~Merca,
Rank partition functions and truncated theta identities,
\textit{Appl. Anal. Discrete Math.} (2021) https://doi.org/10.2298/AADM190401023M.

\bibitem{MK} 
M.~Merca and J.~Katriel,
A general method for proving the non-trivial linear homogeneous partition inequalities, 
\textit{Ramanujan J} \textbf{51}(2) (2020)  245--266. 

\bibitem{MWY}
M.~Merca, C.~Wang and A.~J.~Yee,
A truncated theta identity of Gauss and overpartitions into odd parts, 
\textit{Ann. Comb.}, \textbf{23} (2019) 907--915.

\bibitem{MY}
M.~Merca and A.~J.~Yee,
On the sum of parts with multiplicity at least $2$ in all the partitions of $n$,
\textit{Int. J. Number Theory} \textbf{17}(3) (2021) 665--681.





\bibitem{Penn}
D.~Penniston, 
The $p^a$-regular partition function modulo $p^j$,
\textit{J. Number Theory} \textbf{94} (2002) 320--325.

\bibitem{Penn8}
D.~Penniston,
Arithmetic of $\ell$-regular partition functions,
\textit{Int. J. Number Theory} \textbf{4} (2008) 295--302.





\bibitem{Sloane} 
N.~J.~A.~Sloane, 
The on-line encyclopedia of integer sequences. 
Published electronically at \url{http://oeis.org} (2022).	

\bibitem{Smoot}
N.~A.~Smoot,
On the computation of identities relating partition numbers in arithmetic progressions with eta quotients: an implementation of Radu's algorithm,
\textit{J. Symb. Comput.} \textbf{104} (2021) 276--311.

\bibitem{Subbarao} 
M.~V.~Subbarao, M.~Vidyasagar, 
On Watson's quintuple product identity,
\textit{Proc. Amer. Math. Soc.} \textbf{26} (1970) 23--27.

\bibitem{Xia15}
E.~X.~W.~Xia,
Congruences for some $\ell$-regular partitions modulo $l$,
\textit{J. Number Theory} \textbf{152} (2015) 105--117.

\bibitem{Xia14}
E.~X.~W.~Xia, O.~X.~M. Yao, 
Parity results for $9$-regular partitions, 
\textit{Ramanujan J} \textbf{34} (2014) 109--117.

\bibitem{Xia}
E.~X.~W.~Xia, X.~Zhao,
Truncated sums for the partition function and a problem of Merca,
\textit{Rev. R. Acad. Cienc. Exactas F\'{i}s. Nat. Ser. A Math. RACSAM} \textbf{116}:22 (2022).

\bibitem{Wang17}
L.~Wang,
Congruences for $5$-regular partitions modulo powers of $5$,
\textit{Ramanujan J} \textbf{44}  (2017) 343--358.

\bibitem{Wang}
L.~Wang,
Arithmetic properties of $7$-regular partitions,
\textit{Ramanujan J} \textbf{47} (2018) 99--115.

\bibitem{Webb}
J.~J.~Webb, 
Arithmetic of the $13$-regular partition function modulo $3$,
\textit{Ramanujan J} \textbf{25} (2011) 49--56.






\end{thebibliography}
\end{document}